# ESTIMATION OF THE DENSITY OF REGRESSION ERRORS

By Sam Efromovich[1]

*University of New Mexico*

*Dedicated to the memory of Mark Pinsker, 1925–2003*

Estimation of the density of regression errors is a fundamental issue in regression analysis and it is typically explored via a parametric approach. This article uses a nonparametric approach with the mean integrated squared error (MISE) criterion. It solves a long-standing problem, formulated two decades ago by Mark Pinsker, about estimation of a nonparametric error density in a nonparametric regression setting with the accuracy of an oracle that knows the underlying regression errors. The solution implies that, under a mild assumption on the differentiability of the design density and regression function, the MISE of a data-driven error density estimator attains minimax rates and sharp constants known for the case of directly observed regression errors. The result holds for error densities with finite and infinite supports. Some extensions of this result for more general heteroscedastic models with possibly dependent errors and predictors are also obtained; in the latter case the marginal error density is estimated. In all considered cases a blockwise-shrinking Efromovich–Pinsker density estimate, based on plugged-in residuals, is used. The obtained results imply a theoretical justification of a customary practice in applied regression analysis to consider residuals as proxies for underlying regression errors. Numerical and real examples are presented and discussed, and the S-PLUS software is available.

**1. Introduction.** A residual analysis is a standard part of any regression analysis, and it involves estimation and/or testing of a regression error distribution. This article is devoted to the error density estimation. Let us present the problem, its motivation and possible applications via a classical homoscedastic model, and then more complicated models will be introduced.

Received February 2004; revised September 2004.

[1]Supported in part by BIFAR and by NSF Grant DMS-02-43606.

Following Fan and Gijbels [22], Hart [28] and Eubank [21], suppose that the statistician observes $n$ independent and identically distributed (i.i.d.) realizations $(X_1, Y_1), \ldots, (X_n, Y_n)$ of the pair $(X, Y)$ of random variables. Then the regression problem is to find an underlying regression function $m(x) := E(Y|X=x)$ under the assumption that

$$(1.1) \qquad Y_l = m(X_l) + \xi_l, \qquad l = 1, \ldots, n,$$

$X_1, \ldots, X_n$ are i.i.d. predictors that are uniformly distributed on $[0, 1]$, and $\xi_1, \ldots, \xi_n$ are i.i.d. regression errors that are also independent of the corresponding predictors $X_1, \ldots, X_n$. The model (1.1) is called a homoscedastic regression model with regression errors which are i.i.d. and independent of the predictors. If $\tilde{m}(x)$ is a regression estimate, then $R_l := Y_l - \tilde{m}(X_l)$, $l = 1, \ldots, n$, are called residuals. Patterns in the residuals are used to validate or reject an assumed model. If the model (1.1) is validated, then the next classical step is to look at the distribution of the regression error $\xi$. Because realizations $\xi_1, \ldots, \xi_n$ of regression errors are unavailable to the statistician, residuals are traditionally utilized as their proxies. They may be used either for testing a hypothesis about the underlying error distributions or for estimation/visualization of the error density; see a discussion in the classical text by Neter, Kutner, Nachtsheim and Wasserman [32].

Surprisingly, despite the widespread use of residuals as proxies for unobserved errors, to the best of the author's knowledge, no result about optimal (in any sense) estimation of a nonparametric error density is known. For parametric settings, there exists a recently created Bayesian theory of estimation, and for nonparametric settings, a theory of consistent estimation is developed; the interested reader can find a discussion and further references in [8] and [27]. At the same time, there exists a vast literature devoted to density estimation based on direct observations and to estimation of functionals of the error density; see [2, 14, 34, 37] and [31] where further references can be found.

It is not difficult to understand why the literature on nonparametric error density estimation is practically next to none: the problem is extremely complicated due to its indirect nature. In a nonparametric setting, the difference between any regression estimate and an underlying regression function contains a random term and a bias. The bad news is that additive measurement errors may dramatically slow down optimal rates of density estimation; see [13, 14]. The good news is that, of course, additive errors in residuals become smaller as the sample size increases, and, thus, optimal rates may be preserved. This article shows that, fortunately for applied statistics, the good news prevails under the customary assumption that the regression function is differentiable and the error density is twice differentiable.

It is well known in the nonparametric density estimation literature that rates alone are of little interest for practically important cases of small

DENSITY OF REGRESSION ERRORS 3

datasets, and that rates should be studied together with constants; see the discussion in [30] and [14, 15]. Also, superefficiency and mimicking of oracles are important issues; see the discussion in [4, 5] and [14], Chapter 7. To explore all these issues simultaneously, it is convenient to employ an oracle approach suggested by Mark Pinsker more than two decades ago. Namely, suppose that an oracle (which will be referred to as a Pinsker oracle and a particular one is defined in Appendix B) knows the underlying regression errors $\{\xi_l, l = 1, \ldots, n\}$ and the oracle possesses a bouquet of desired statistical properties like sharp minimaxity, superefficiency, matching more powerful oracles that know an estimated error density, and so on. Then, if the statistician can suggest a data-driven error density estimate that matches the Pinsker oracle, this estimator simultaneously solves all the above-formulated problems. Moreover, Pinsker conjectured that a plug-in Pinsker oracle, based on residuals, may be the wished data-driven estimator. This article proves this long-standing Pinsker conjecture and, as a particular corollary, establishes minimax rates and constants of the error density estimation.

There are many practical applications of the error density estimation. Let us mention a few that will guide us in this article. (i) Well-known classical applications are data interpretation, inference, decision making, hypothesis testing, the diagnostics of residuals, model validation and, if necessary, model adjustment in terms of the error distribution. (ii) Another classical application, which is actually the main aim of any regression analysis, is the prediction of a new observation where the error density plays the pivotal role; see [32], Section 2.5. (iii) Goodness-of-fit tests are another natural application; see [2] and [28]. (iv) The error density is used in a sharp minimax regression estimation; see [12]. (v) The error density can be used in statistical quality control and classification; see [16], as well as a discussion in Section 2.

The model (1.1) with a uniformly distributed predictor is definitely the most frequently studied in the regression literature, but a regression analysis may reveal patterns that contradict this simple model. For instance, predictors may not be uniform and/or the errors may have different variances. In this case either some remedial measures like a data transformation and/or weighting are applied (these remedies are not discussed here and the interested reader is referred to the books by Carroll and Ruppert [6] or Neter, Kutner, Nachtsheim and Wasserman [32]), or model (1.1) with an unknown design density $p(x)$ is considered, or a more general heteroscedastic regression model is considered:

$$(1.2) \qquad Y_l = m(X_l) + \sigma(X_l)\xi_l, \qquad l = 1, \ldots, n,$$

where $\sigma(x)$ is a (positive) scale function, the errors $\{\xi_1, \ldots, \xi_n\}$ are i.i.d. with zero mean, unit variance and independent of the corresponding predictors, and the predictors are i.i.d. according to an unknown design density $p(x)$ supported on $[0,1]$. Following Pinsker's paradigm, a data-driven error density estimator should be compared with an oracle that knows the underlying errors $\xi_1, \ldots, \xi_n$. Here it is natural to use rescaled residuals as proxies for unobserved errors, and an implementation of this path implies estimation of both regression and scale functions and then dealing with additive and multiplicative measurement errors. It will be shown that, under the assumption of the differentiability of each nuisance function and a known finite support of an estimated error density, the plug-in Pinsker oracle still matches the Pinsker oracle; the case of errors with infinite support is an open problem.

Now we are in position to consider another assumption about models (1.1)–(1.2) that needs to be addressed: independence between regression errors and predictors. There are many known examples where this assumption does not hold; see particular ones in Section 2. Another customary situation is where an underlying model is heteroscedastic, but the statistician assumes/believes that it is homoscedastic; an interesting particular example is presented in [28], pages 257–258. If we simply ignore a possible dependence between $X$ and $\xi$, then what does our plug-in estimate exhibit or, in other words, what do residuals proxy? To the best knowledge of the author, there is no nonparametric literature devoted to this issue. This article shows that in this case the marginal error density is estimated and then all the above-discussed statistical results hold. In particular, this establishes that a plug-in estimation is robust toward a possible dependency between predictor and regression error, and this is an important conclusion for an applied residual analysis.

Finally, let us note that the developed theory of plug-in estimation significantly simplifies the problem of creating software because known statistical programs can be used directly. This article uses the S-PLUS software package of Efromovich [14] which is available on request from the author.

The structure of this article is as follows. Section 2 presents several numerical simulations and real practical examples that should help the reader to understand the problem, its solution and possible applications. Section 3 contains mathematical results, and discussion is presented in Section 4. Appendix A describes the main steps of proofs; complete proofs can be found in [16, 17, 19]. Appendix B is devoted to the Pinsker oracle, and it presents new results for the case of densities with infinite support.

**2. Several examples.** Let us explain the above-described problem of error density estimation via several particular examples.

Figure 1 presents a simulation conducted according to model (1.2) with functions described in the caption. The left-top diagram exhibits a scattergram, and the problem is to estimate an underlying error density. Asymptotic theory, presented in Section 3, indicates that the S-PLUS software package of Efromovich [14] can be used for calculating rescaled residuals and then error density estimation. Recall that the package supports Efromovich–Pinsker (EP) adaptive series estimates; see the discussion of regression, scale and density estimates in [14], Sections 4.2, 4.3 and 3.1. Let us explain how this software performs for the simulated dataset. The scattergram is overlaid by the EP regression estimate (the dotted line) and it can be compared with the underlying regression $m(x)$ (the solid line). This particular estimate is not perfect and we can expect relatively large additive measurable errors in the residuals. The left-middle diagram exhibits the EP scale estimate (the dotted line), and it can be compared with the underlying scale function $\sigma(x)$ (the solid line). This estimate is also not perfect, so we can expect multiplicative measurement errors in the rescaled residuals shown in the left-bottom diagram. The right column of diagrams exhibits the process of the error density estimation by the Pinsker oracle and the corresponding plug-in estimation. The Pinsker oracle is based on unobserved errors shown in the right-top diagram, and the plug-in estimate is based on rescaled residuals shown in the right-middle diagram. The oracle, the plug-in estimate and the underlying error density are shown in the right-bottom diagram by dashed, dotted and solid lines, respectively.

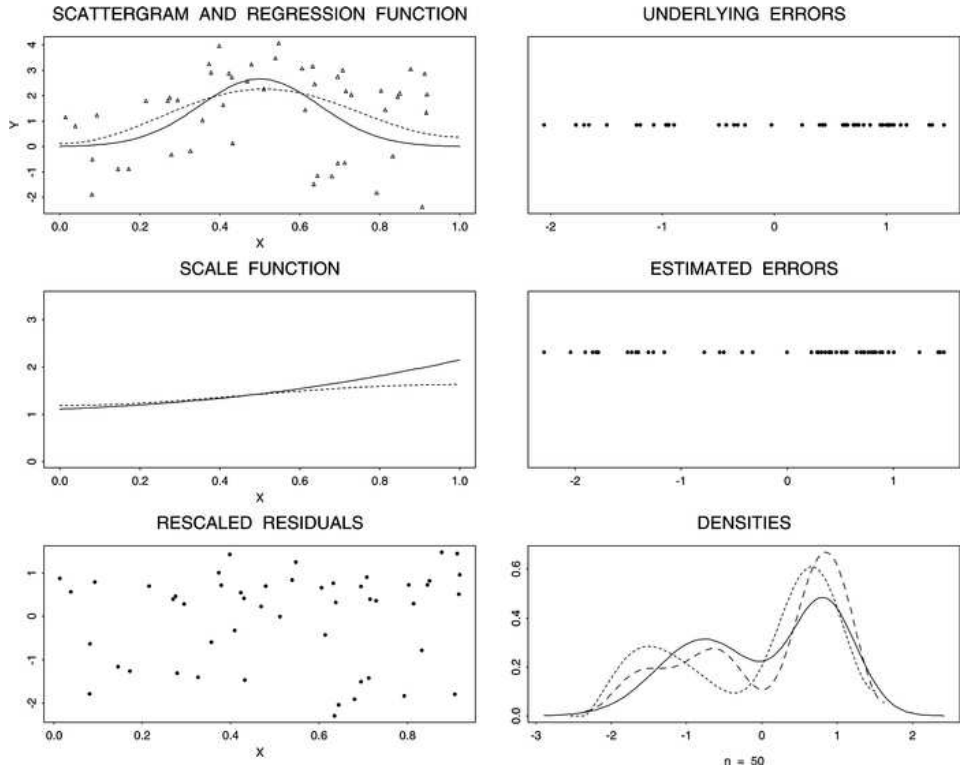

Fig. 1. *Simulated example of heteroscedastic regression* (1.2) *with the regression function being the Normal, the scale function being the Monotone, the design density being the Uniform and the error density being the Bimodal; these underlying functions are defined in* [14], *page* 18. *The simulated scattergram is shown by triangles, the sample size* $n = 50$ *is shown in the subtitle of the right-bottom diagram. The dotted lines show data-driven estimates, the solid lines show underlying functions, and the dashed line in the right-bottom diagram shows the oracle estimate based on underlying errors exhibited in the right-top diagram.*

As we see, due to the presence of measurement errors, the data-driven estimate performs worse than the oracle. The estimate is overspread, and this outcome is typical for data contaminated by measurement errors; see the discussion in [13] and [14], Section 3.5. Nonetheless, the estimate correctly indicates the bimodal nature of the error. Keeping in mind that any nonparametric analysis is considered as a first glance at the data, the estimate readily indicates that the error density is not normal. This conclusion implies that classical methods of regression analysis, based on normal distribution of errors, should be modified. For instance, a prediction error may be described by using the error density estimate.

Let us complement this single simulation with an intensive Monte Carlo study where 500 identical simulations are conducted for each sample size from the set $\{25, 50, 75, 100, 150, 200\}$. For each simulation, we calculate the ratio of ISEs of the Pinsker oracle and the estimate, and then for 500 simulations, calculate the sample mean, sample median and sample standard deviation of the ratios. The corresponding results are as follows: $\{(1.05/0.93/0.74); (1.01/0.83/0.72); (0.96/0.81/0.64); (0.97/0.85/0.63); (0.94/0.88/0.53); (0.96/0.87/0.56)\}$, where an element $(A/B/C)$ presents the sample mean, median and standard deviation for a corresponding sample size. Note that a mean ratio or median ratio smaller than 1 favors the Pinsker oracle. As we see, for the explored sample sizes, traditionally considered as small even for the case of direct observations, plug-in estimation performs respectively well. This tells us that Pinsker's proposal of comparing a data-driven estimator with an oracle is feasible even for the smallest samples. The interested reader can find more simulations and numerical studies in [16].

Our next simulation, exhibited in Figure 2, addresses an important issue of rescaling of residuals. It is fair to say that an applied regression analysis is primarily devoted to a homoscedastic regression, and a possible issue of heteroscedasticity is addressed by a data transformation and/or weighting rather than rescaling; see the discussions in [6] and [32]. We shall comment on this shortly, but now let us consider an example of a homoscedastic regression (1.1) which is treated by the suggested software that always attempts to rescale residuals. A simulated scattergram is shown in the left-top diagram of Figure 2. Absolute values of residuals are shown by squares in the left-middle diagram, and they readily exhibit heteroscedasticity. We know that this heteroscedasticity is stochastic in nature (look at the underlying horizontal scale function shown by the solid line), but the software does not know this. Thus, it is of interest to understand how the software will perform with respect to a new oracle that knows that the model is homoscedastic. In other words, let us compare performances of the same density estimator where scaled residuals and residuals are plugged in. The left-bottom and right-middle diagrams exhibit by dots and squares rescaled residuals and residuals, respectively. The corresponding density estimates are shown by the dotted and long-dashed lines in the right-bottom diagram; the solid and dashed lines in this diagram exhibit the underlying normal error density and the Pinsker oracle (based on unobserved regression errors), respectively. As we see, in this particular case the rescaling had a positive effect on the estimation. In general, this cannot be the case, so let us conduct a numerical study identical to the above-described one with the only difference being that now we are comparing the use of rescaled residuals (the estimate) and residuals (a new oracle). The results are the following: $\{(0.99/0.91/0.78); (0.97/0.78/0.57); (1.02/0.78/0.65); (0.93/0.87/0.61); (0.98/0.88/0.57); (0.97/0.87/0.54)\}$. The study indicates that rescaling can be considered a robust procedure for homoscedastic regression, and Section 3 presents asymptotic justification of this empirical observation.

The main purpose of the next simulation is to allow us to discuss the case of error depending on the predictor, and it also allows us to explore possible applications for statistical quality control. Assume that a process is inspected at ordered times $X_l$ and corresponding observations are $Y_l$, $l = 1, \ldots, n$. Recall that it is customary to say that a process is in control if its mean (centerline, regression function) and standard deviation (scale function, volatility) are constant. Keeping in mind that a traditionally assumed distribution of controlled variables is Gaussian, the latter implies a stationary distribution of the process; see the discussion in [10], Chapter 23. The two top diagrams in Figure 3 present a simulated process together with its two main characteristics. Because mean and standard deviation of the process are constant, the process is declared to be in control. However, even if the process is in control, it may not be strictly stationary. Thus, let us continue our analysis of the process. The third diagram shows us the estimated marginal density of residuals (the dotted line), which exhibits a non-Gaussian shape (note that the underlying marginal density is shown by the solid line). If it is known that the process must be Gaussian, this error density raises a red flag. If no action is required, as in the familiar "normal tool wear" example, then modified acceptance charts and hypotheses tests, based on the estimated density, should be suggested; see [10], Chapters 23 and 25. To check the drawn conclusion about the changed distribution, the two bottom diagrams exhibit an onset error density for the first 50 observations and an sunset error density for the last 50 observations. They support our preliminary conclusion that the error distribution is changing. This example shows that nonparametric error density analysis can be a valuable addition to classical quality control methods.

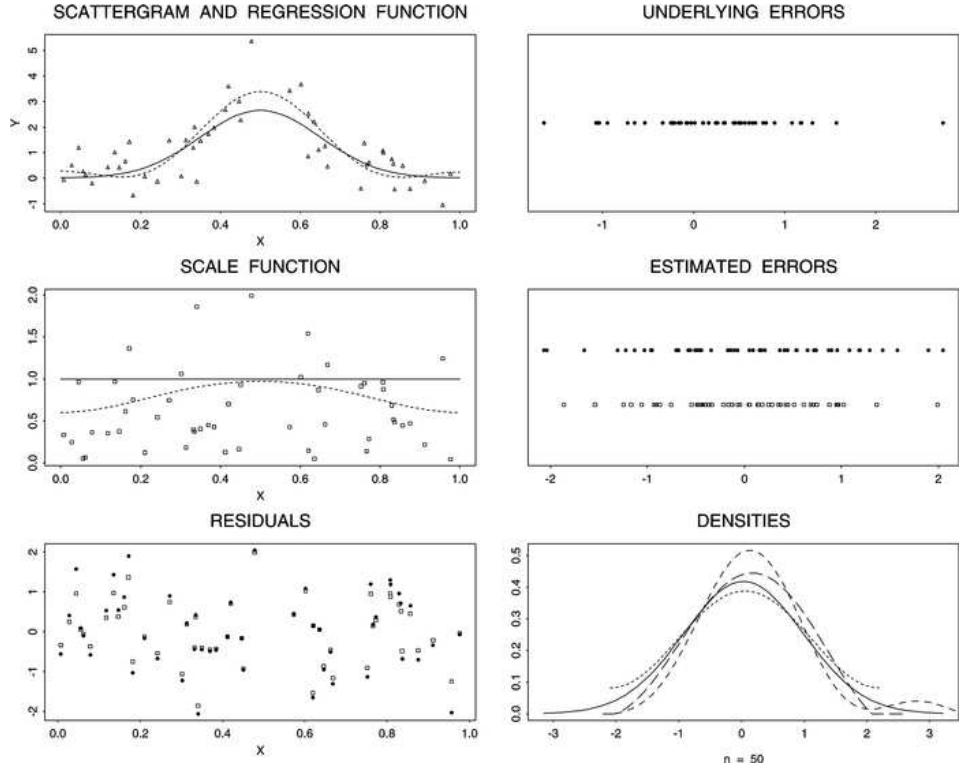

FIG. 2. *The use of rescaled residuals and residuals in a homoscedastic regression. The structure of the diagrams is similar to Figure* 1 *with the following modification. Rectangles in the left-middle diagram show absolute values of residuals. Rectangles and dots in the left-bottom diagram and the right-middle diagram exhibit residuals and rescaled residuals, respectively. The long-dashed line in the right-bottom diagram exhibits the estimate based on residuals.*

Now we are in a position to explore several real practical examples. The research division of BIFAR, a company with businesses in equipment and chemicals for wastewater treatment plants, has studied performance of a centrifuge for mechanical dewatering of a sludge produced by a food processing plant. The aim of the study has been to understand how a sludge, containing a fat waste, can be centrifuged. The top-left diagram in Figure 4 presents results of a particular experiment. Index of fat is the predictor and index of centrifuging is the response. It has been established in [16] that the distribution of regression errors crucially depends on the predictor. Thus, we know a priori that we will visualize the marginal error density.

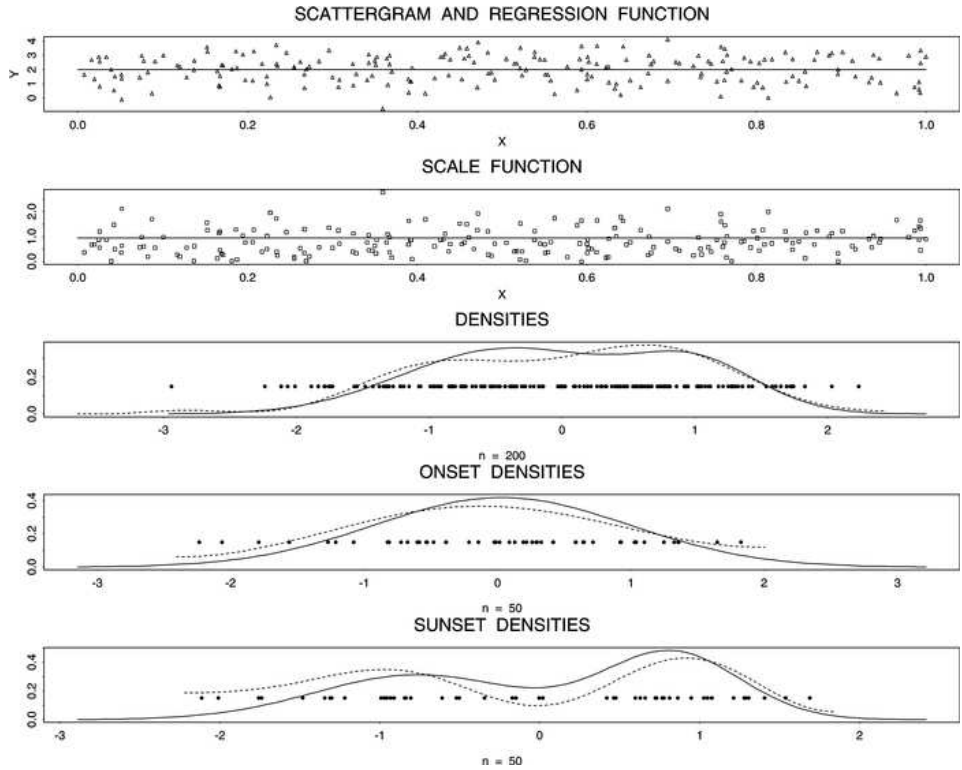

FIG. 3. *Simulated example with the error distribution depending on the predictor. Here $Y = 2 + \varepsilon(X)$ where the error is a linear mixture of the Bimodal density, shown in Figure 1 and having weight $X$, and the Normal density, shown in Figure 2 and having weight $1 - X$. The structure of the two top diagrams is similar to the left ones in Figure 2. The third diagram exhibits the estimated and underlying marginal densities. The two bottom diagrams show the marginal error density estimates for initial (onset) and final (sunset) 50 observations. The estimates and underlying densities are shown by dotted and solid lines, respectively.*

Before discussion of the example, let us make the following remark about the software. It allows the statistician to estimate error densities with a known manually chosen finite support or infinite/unknown support. Intensive simulations in [16] show that, for smaller sample sizes, the former approach benefits the estimation, while, for larger samples, both methods perform similarly. In the simulated examples support has been unknown and, thus, the shown estimates are completely data-driven. A manual choice of support is not a difficult step in many applied settings because it is defined by specifications. In particular, for the BIFAR example, this approach implied the manual choice $[-2.75, 2.75]$ for the support. Due to the small sample size $n = 47$, this help is valuable and should be utilized whenever possible. (The interested reader can find discussion of several manual and data-driven choices of support in [14], Chapter 3.)

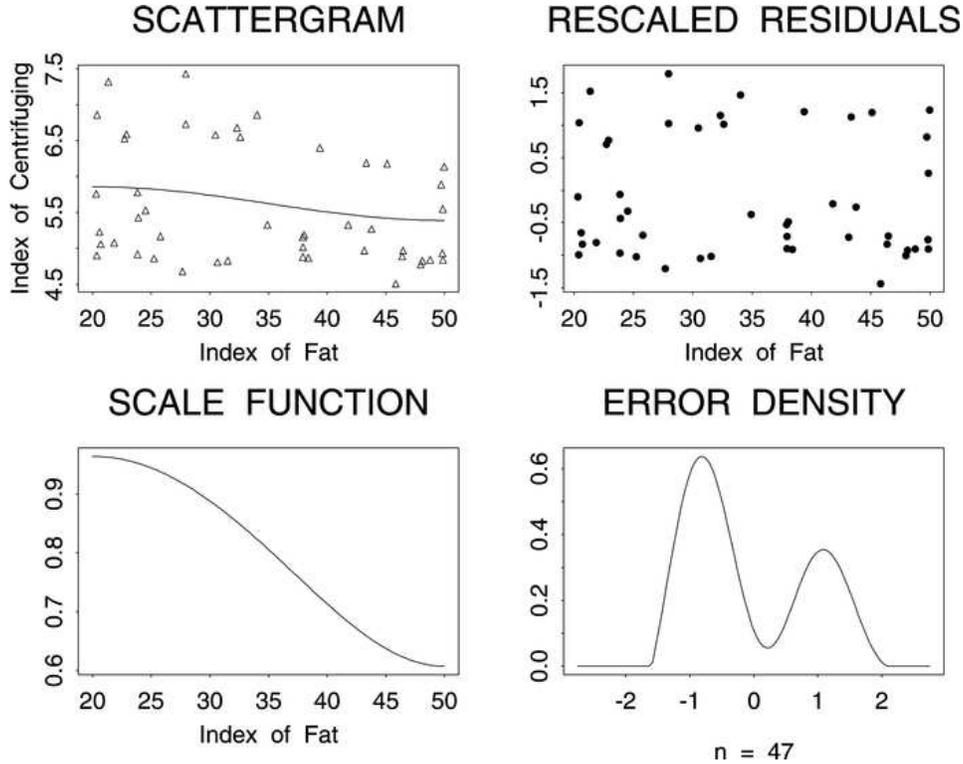

FIG. 4. *Centrifuging a food processing plant's waste. The sample size is $n = 47$. The scattergram is shown by triangles and it is overlaid by the EP regression estimate.*

Now let us explore the BIFAR dataset. The left-top diagram in Figure 4 exhibits the data and the estimated regression function. It is readily seen from this diagram that the regression is heteroscedastic. The bottom-left diagram contains the scale estimate, and it supports our visual conclusion about the heteroscedasticity. Let us note that neither the regression nor the scale estimate has been a surprise for BIFAR. The right-top diagram shows rescaled residuals; the diagram indicates that regression and scale estimates performed well and no heteroscedasticity can be observed. Also, after a closer look at the rescaled residuals, it is possible to note clusters in the residuals. This observation is supported by the estimated marginal density of errors shown in the right-bottom diagram. The density estimate reveals that it is a mixture of two distributions with the larger "left" cluster having a negative bias which "drags" the index of centrifuging down. This was a fantastic insight into the centrifuging process for BIFAR that, just for free, gave the company a new tool for the process analysis. As a result, while classical regression analysis traditionally describes a relationship between two variables by univariate regression and scale functions, it is proposed to complement the analysis by an extra univariate function—error density. Let us stress that it would be great to complement this analysis with a conditional density, but the sample is too small for bivariate function estimation.

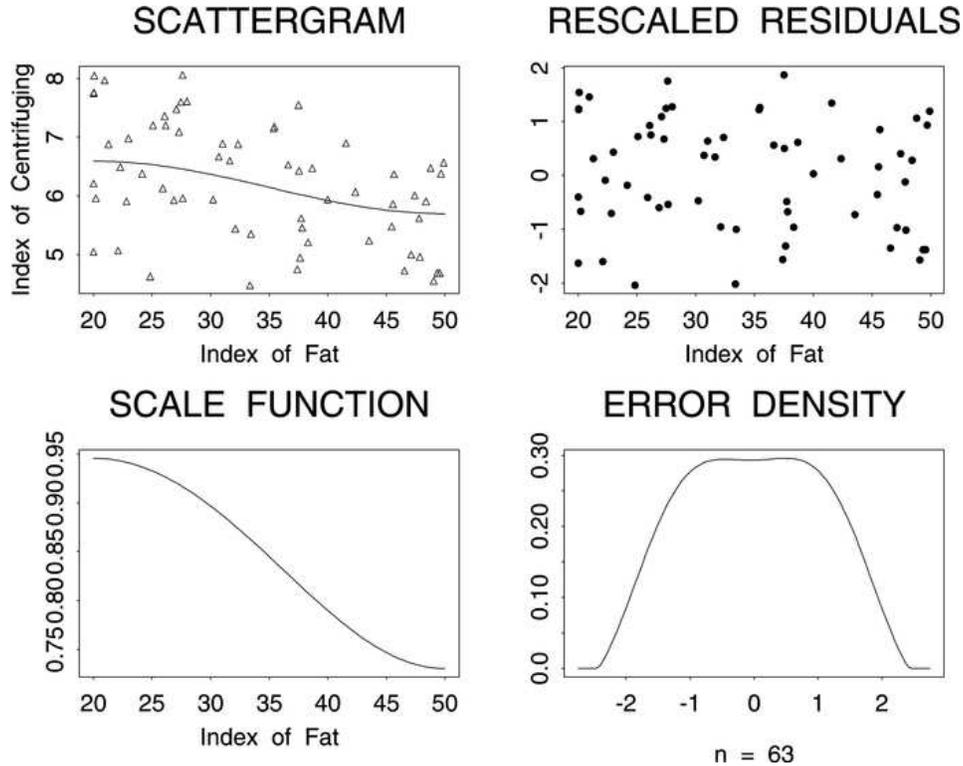

Fig. 5. *Centrifuging a food processing plant's waste with added chemical. Sample size $n = 63$.*

Based on this outcome, BIFAR decided to conduct a series of experiments where special chemicals were added to the sludge. Figure 5 presents (in the same format) results of a particular experiment. Note that the regression and scale functions are about the same due to robust performance of the centrifuge. On the other hand, the marginal error density indicates that the chemical is able to merge together the above-discussed clusters, and it also decreases the relative effect of the "left" cluster. This may explain, at least partially, the overall increase in the index of centrifuging caused by the chemical. Note that all these observations have been based on the analysis of univariate functions. Of course, it would be nice to evaluate an underlying conditional density, but the dataset is too small for this and we are restricted to the univariate nonparametric analysis.

Let us make a comment that connects these two practical examples with the simulated example in Figure 3. It is possible to imagine a situation where BIFAR observations are a time series where, due to some circumstances, index of fat increases. For instance, this may occur if a food processing plant illegally dumps its waste into a municipal sewage system. Then the error density will be the first indicator of such violation. Also, these two figures indicate the possibility of using the error density for classification purposes: the chemical is present or not, fat is present or not. In other words, an error density is an additional univariate characteristic (in addition to mean and scale functions) that may be useful in many settings of industrial statistics.

We may conclude that, on the top of such classical applications in regression analysis as prediction, model validation, hypothesis testing and optimal estimation of regression functions, error-density estimation is a valuable and feasible data-analysis tool on its own in time series, quality control and industrial statistics.

**3. Optimal estimation of the error density.** The aim of this section is twofold. First, we would like to establish minimax rates and constants of mean integrated squared error (MISE) convergence of error density estimates in homoscedastic models (1.1) and, if possible, in heteroscedastic models (1.2) with errors depending on predictors. Recall that, even for a homoscedastic model, minimax rates are unknown; see [8]. Second, we would like to suggest a data-driven (adaptive) estimator that attains the minimax convergence. Ideally, to support the classical methodology of applied regression analysis and to employ available statistical software, such an estimator should be a known (for direct observations) density estimator based on appropriately calculated residuals, that is, it should be a plug-in density estimator.

Two classical models of errors will be studied: models of errors with a known finite support $[a, a+b]$ and errors with infinite support $(-\infty, \infty)$. Recall that we discussed particular examples in Section 2. We need to make a comment about the finite interval case. It will be convenient to evaluate the density over a fixed interval, and a customary interval is $[0,1]$. In models (1.1) and (1.2) the error support cannot be $[0,1]$ because $E\xi = 0$. Thus, we employ a familiar location-scale transformation, introduce a new random variable $\varepsilon := (\xi - a)/b$ and then study the equivalent problem of estimation of the density $f^\varepsilon(u)$ of the transformed error $\varepsilon$ instead of the density $b^{-1}f^\varepsilon([u-a]/b)$ of $\xi$. The approach of estimation of a rescaled random variable is discussed in detail in [14], Chapter 3. From now on we omit the superscript $\varepsilon$ in the density, denote by $f$ the density of $\varepsilon$, and refer to it as the error density (of interest).

In what follows, with some obuse of notation, we shall always present results for finite and infinite supports simultaneously.

3.1. *Model and assumption for the case of finite support.* The studied regression model is (1.2), where the error $\xi$ may depend on the predictor $X$. Neither the regression function $m$ nor the scale function $\sigma$ is supposed to be known. The observed predictors $(X_1, \ldots, X_n)$ are i.i.d. according to an unknown design density $p(x)$ supported on $[0,1]$; the regression errors $\xi_1, \ldots, \xi_n$ do not take on values beyond a known finite interval $[a, a+b]$ and may depend on the corresponding predictors according to an unknown conditional density $b^{-1}\psi([\nu - a]/b|x)$, $\nu \in [a, a+b]$; the pairs $(X_1, \xi_1), \ldots, (X_n, \xi_n)$ are supposed to be independent and identically distributed. The marginal density of the rescaled errors $\varepsilon_l = [\xi_l - a]/b$ is the object of interest; that is, the issue is to estimate the density $f(u) = \int_0^1 \psi(u|x) p(x)\,dx$, $u \in [0,1]$, based on $n$ pairs of observations $\{(X_1, Y_1), \ldots, (X_n, Y_n)\}$.

ASSUMPTION A. The regression function $m(x)$, the design density $p(x)$ and the scale function $\sigma(x)$ are differentiable and their derivatives are bounded and integrable on $[0,1]$. Also, $\min_{x \in [0,1]} \min(\sigma(x), p(x)) > 0$ and $\int_0^1 p(x)\,dx = 1$.

ASSUMPTION B (Finite support). Model (1.2) is considered where the errors may depend on the predictors. Pairs of observations $(X_1, Y_1), \ldots, (X_n, Y_n)$ are i.i.d. The conditional density $\psi(u|x)$ is such that $\frac{\partial}{\partial x} \frac{\partial^2}{\partial u^2} \psi(u|x)$ exists, is bounded and integrable on $[0,1]^2$, and $\psi(u|x) = 0$ for $u \notin (0,1)$, $x \in [0,1]$.

ASSUMPTION C. For i.i.d. observations (errors) $Z_1, \ldots, Z_r$ from a density $f(u)$ with unit support $[0,1]$ or infinite support $(-\infty, \infty)$, Appendix B defines a data-driven density estimate $\hat{f}_P(u, Z_1^r)$, $Z_1^r := (Z_1, \ldots, Z_r)$. This estimate, based on underlying errors, is employed as the Pinsker oracle. It is assumed that the statistician knows all parameters of this estimate.

3.2. *Model and assumption for the case of infinite support.* Due to the complexity of the case, the studied model is homoscedastic regression (1.1) where the error $\xi$ is independent of the predictor $X$. Neither the regression function $m$ nor the design density $p$ of the predictors is known. The problem is to estimate the density $f(u)$ of $\xi$ based on $n$ i.i.d. pairs of observations $(X_1, Y_1), \ldots, (X_n, Y_n)$. In what follows a reference to the above-formulated Assumption A means that $\sigma(x) = 1$, $x \in [0,1]$.

ASSUMPTION B (Infinite support). Model (1.1) is considered with the error being independent of the predictor and pairs of observations $(X_1, Y_1), \ldots,$ $(X_n, Y_n)$ being i.i.d. The error density $f(u)$ is supposed to be square integrable, that is, $\int_{-\infty}^{\infty} f^2(u)\,du < \infty$, and its characteristic function $h(v) := \int_{-\infty}^{\infty} f(u)e^{iuv}\,du$ satisfies $\int_{-\infty}^{\infty} v^4|h(v)|^2\,dv < \infty$.

3.3. *Notational convention.* Several sequences in $n$ are used: $b_n := 4 + \ln\ln(n+20)$; $n_2 := n - 3n_1$; $n_1$ is the smallest integer larger than $n/b_n$; $S := S_n$ is the smallest integer larger than $n^{1/3}$. In what follows we always consider sufficiently large $n$ such that $\min(n_1, n_2) > 4$. $C$'s denote generic positive constants, $o(1) \to 0$ as $n \to \infty$, and integrals are taken over $[0, 1]$ or $(-\infty, \infty)$, depending on the support considered. Also, $(x)_+ := \max(0, x)$.

3.4. *Plugged-in residuals.* The aim of this section is to explain a procedure for the calculation of plugged-in residuals. Four different subsamples are used to estimate the design density, the regression function, the scale function and the error density, respectively (the author conjectures that all $n$ observations may be used for estimation of each function and the result will still hold). The first $n_1$ observations are used to estimate the design density $p(x)$, the next $n_1$ observations are used to estimate the regression function $m(x)$, the next $n_1$ observations are used to estimate the scale function $\sigma(x)$, and the last $n_2$ observations are used to estimate the error density of interest $f(u)$. Note that $n_2 \geq [1 - 3(b_n^{-1} + n^{-1})]n$ and, thus, using either $n_2$ or $n$ observations implies the same MISE convergence. The design density estimate $\hat{p}(x)$ is a truncated cosine series estimate,

$$\hat{p}(x) = \max\left(b_n^{-1}, n_1^{-1} \sum_{s=0}^{S} \sum_{l=1}^{n_1} \varphi_s(X_l)\varphi_s(x)\right).$$

The regression estimate $\hat{m}(x)$ is also a truncated cosine series estimate,

$$\hat{m}(x) = \sum_{s=0}^{S} \hat{\kappa}_s \varphi_s(x), \tag{3.1}$$

where

$$\hat{\kappa}_s = n_1^{-1} \sum_{l=n_1+1}^{2n_1} Y_l \hat{p}^{-1}(X_l) \varphi_s(X_l). \tag{3.2}$$

Under model (1.2), the scale estimate $\hat{\sigma}(x)$ is also a truncated cosine series estimate,

$$\hat{\sigma}(x) = \min(\max(\tilde{\sigma}(x), b_n^{-1}), b_n), \tag{3.3}$$

where $\tilde{\sigma}(x) = \sqrt{(\tilde{\sigma}^2(x))_+}$ and $\tilde{\sigma}^2(x)$ is a regression estimate defined identically to (3.1)–(3.2), where pairs $\{(X_l, Y_l), l = n_1+1, \ldots, 2n_1\}$ are replaced by $\{(X_l, [Y_l - \hat{m}(X_l)]^2), l = 2n_1+1, \ldots, 3n_1\}$.

Then, for finite support (recall that in this case a heteroscedastic model is considered) we define rescaled residuals

$$\hat{\varepsilon}_l := \frac{Y_l - \hat{m}(X_l)}{b\hat{\sigma}(X_l)} - \frac{a}{b}, \qquad l = n - n_2 + 1, \ldots, n. \tag{3.4}$$

For infinite support (in this case a homoscedastic model is considered) we define residuals

$$\hat{\xi}_l := Y_l - \hat{m}(X_l), \qquad l = n - 2n_1 + 1, \ldots, n. \tag{3.5}$$

Now we can use a unified notation for the residuals and underlying errors. Denote by $\hat{\mathbf{Z}}$ a vector $(\hat{\varepsilon}_{n-n_2+1}, \ldots, \hat{\varepsilon}_n)$ or a vector $(\hat{\xi}_{n-2n_1+1}, \ldots, \hat{\xi}_n)$ for finite and infinite support cases, respectively. Similarly, $\mathbf{Z}$ denotes a vector of transformed errors $(\varepsilon_1, \ldots, \varepsilon_n)$ or a vector of errors $(\xi_1, \ldots, \xi_n)$ for finite and infinite support cases, respectively. Note that $\mathbf{Z}$ is known to the Pinsker oracle but not to the statistician.

3.5. *Main assertion.* It is possible to show that, under the given assumptions, the MISE of the plug-in Pinsker oracle $\hat{f}_P(u, \hat{\mathbf{Z}})$, defined in Appendix B, can asymptotically match the MISE of the Pinsker oracle $\hat{f}_P(u, \mathbf{Z})$ based on underlying regression errors.

THEOREM 1. *The cases of finite and infinite supports are considered simultaneously. Suppose that Assumptions A, B and C hold. Then, for all sufficiently large samples such that* $\min(n_1, n_2) > 4$, *the MISE of the plug-in Pinsker oracle satisfies the Pinsker oracle inequality*

$$\begin{aligned} E \int (\hat{f}_P(u, \hat{\mathbf{Z}}) - f(u))^2 \, du \\ \leq (1 + P^* \ln^{-1}(b_n)) E \int (\hat{f}_P(u, \mathbf{Z}) - f(u))^2 \, du + P^* b_n^3 n^{-1}, \end{aligned} \tag{3.6}$$

*where $P^*$ is a finite constant.*

Recall that $b_n = 4 + \ln \ln(n + 20)$ and, thus, $P^* b_n^3 n^{-1} = o(1) \ln(n) n^{-1}$, that is, the second term in (3.6) is negligible with respect to minimax MISEs of analytic and differentiable densities which are at least of order $\ln(n) n^{-1}$. Also note that Assumptions A and B involve no interplay between smoothness of the error density and smoothness of the triplet of nuisance functions (design density, regression and scale). This allows us to conclude that residuals can be considered as proxies for unobserved regression errors, and this conclusion supports the customary methodology of applied statistics.

The obtained result also allows us to establish minimax rates and constants of MISE convergence; they are presented below.

3.6. *Optimal rates and constants of MISE convergence.* This section answers several classical questions about optimal estimation of a nonparametric error density. To the best of the author's knowledge, so far no results about optimal rates have been known even for the simplest homoscedastic regression model (1.1) with uniformly distributed predictors.

Here we are considering a Sobolev ($\alpha$-fold differentiable) class $\mathcal{S}(\alpha, Q)$ of error densities and an analytic class $\mathcal{A}(\gamma, Q)$ of error densities. These classes are defined and discussed in Appendix B for finite and infinite supports, and let us note that the same notation is used in both cases.

COROLLARY 1 (Differentiable error density). *Suppose that the assumptions of Theorem 1 and* (B.12) *of Appendix B hold and $\alpha \geq 2$. Then the plug-in Pinsker oracle is sharp minimax over Sobolev error densities and all possible oracles, that is,*

$$\sup_{f \in \mathcal{S}(\alpha, Q)} E \int [r_n(\mathcal{S}(\alpha, Q))(\hat{f}_P(u, \hat{\mathbf{Z}}) - f(u))]^2 \, du$$

(3.7) $$= (1 + o(1)) \inf_{\check{f}} \sup_{f \in \mathcal{S}(\alpha, Q)} E \int [r_n(\mathcal{S}(\alpha, Q))(\check{f}(u, \mathbf{Z}) - f(u))]^2 \, du$$

$$= (1 + o(1)),$$

*where the infimum is taken over all possible oracles $\check{f}$ based on unavailable-to-the-statistician errors $\mathbf{Z}$ and parameters $\alpha$ and $Q$, the sharp normalizing factor is*

(3.8) $$r_n(\mathcal{S}(\alpha, Q)) := [n^{2\alpha/(2\alpha+1)} / P(\alpha, Q)]^{1/2}$$

*and $P(\alpha, Q)$ is the famous constant of Pinsker* [35],

(3.9) $$P(\alpha, Q) := (2\alpha + 1)[\pi(2\alpha + 1)(\alpha + 1)\alpha^{-1}]^{-2\alpha/(2\alpha+1)} Q^{1/(2\alpha+1)}.$$

COROLLARY 2 (Analytic error density). *Suppose that the assumptions of Theorem 1 and* (B.12) *of Appendix B hold. Then the plug-in Pinsker oracle is sharp minimax over analytic error densities and all possible oracles, that is,*

$$\sup_{f \in \mathcal{A}(\gamma, Q)} E \int [r_n(\mathcal{A}(\gamma, Q))(\hat{f}_P(u, \hat{\mathbf{Z}}) - f(u))]^2 \, du$$

(3.10) $$= (1 + o(1)) \inf_{\check{f}} \sup_{f \in \mathcal{A}(\gamma, Q)} E \int [r_n(\mathcal{A}(\gamma, Q))(\check{f}(u, \mathbf{Z}) - f(u))]^2 \, du$$

$$= (1 + o(1)),$$

*where the infimum is taken over all oracles $\check{f}$ based on unavailable-to-the-statistician errors $\mathbf{Z}$ and parameters $\gamma$ and $Q$, and the sharp normalizing factor is*

(3.11) $$r_n(\mathcal{A}(\gamma, Q)) := [(2\pi\gamma)n / \ln(n)]^{1/2}.$$

The results establish that, whenever Assumptions A and B hold, minimax rates and constants of MISE convergence for the error density estimation are the same as for the case of directly observed errors. Moreover, the minimax estimator is a plug-in one based on appropriately calculated residuals, and it satisfies the oracle inequality. These results verify the long-standing Pinsker conjecture.

## 4. Discussion.

4.1. It is an important fact that Assumption A (about properties of the regression function, scale function and design density) and Assumption B (about properties of an estimated error density) do not interplay. Also, the minimal restrictions on smoothness of all these functions are classical in the nonparametric literature; see [21, 22, 28, 37].

4.2. The assumption $\int_{-\infty}^{\infty} v^4 |h(v)|^2 \, dv < \infty$ about the characteristic function in Assumption B (infinite support) is identical to the assumption that the second generalized derivative of $f(u)$ is square integrable; see [33], page 35. Thus, the assumptions for error densities with finite and infinite supports are similar.

4.3. Let us heuristically explore the presented results from the point of view of the prediction of a new observation $Y^*$ at a random level $X$ of the predictor; see [32], Section 2.5. Whatever prediction topic is considered (hypothesis testing, confidence intervals, etc.), the error density plays a crucial role. Consider the classical model (1.1), and recall that a traditional applied approach/paradigm is to assume that $Y^* = \hat{m}(X) + \eta$, where $\hat{m}$ is a regression estimate, $\eta$ is a prediction error with a density $\tilde{f}$, and the regression and error density estimates are based on the previous $n$ observations. The prediction problem resembles the one considered in the article, so it is natural to explore how the regression and error density estimates suggested in Section 3 will work together in the prediction problem. We note that, according to (1.1), the prediction error can be written as $\eta = m(X) - \hat{m}(X) + \xi$; thus, to verify the paradigm "$\xi$ mimics $\eta$," we need to understand how the difference $m(X) - \hat{m}(X)$ affects the density of $\eta$. Recall that this difference has a classical decomposition into a zero-mean random component and bias. To simplify the heuristic, let us consider only the effect of bias; denote the bias as $b(X)$. Under Assumption A, the squared bias can be (at most) of order $n^{-2/3}$. This implies that the characteristic function $h_b(v) := E\{e^{ivb(X)}\}$ of the bias is close to 1 for frequencies $|v| < o(1)n^{1/3}$, and we can conclude that on these frequencies the characteristic function of $\xi$ does mimic the characteristic function of $\eta$. [Note that beyond these frequencies the characteristic function $h_b(v)$ may be separated from 1.] Recall that at least twice-differentiable densities, considered in this article, are estimated only on frequencies $|v| \leq O(1)n^{1/5}$. As a result, the paradigm holds (of course, we have considered only the bias effect, but similar arguments are applied to the random component). On the other hand, let us relax Assumption B and assume that the error density is only differentiable. Then a rate-optimal error density estimation requires evaluation of its characteristic function on frequencies $|v| \leq O(1)n^{1/3}$, and then the distributions of $\eta$ and $\xi$ may be different. Of course, in this case we can employ a nonoptimal error density estimation, which involves evaluation of the characteristic function only on frequencies $|v| < o(1)n^{1/3}$. The latter preserves the paradigm at the expense of the error density estimation. What we have observed is the onset of irregularity in the error density estimation, and this is an interesting and challenging topic on its own.

4.4.　　There will be a separate paper about the case of infinite support and heteroscedastic regression. Due to the presence of multiplicative measurement errors in residuals, this case requires an additional assumption on the tails of the error distribution, and it is a technically involved matter to suggest a mild assumption.

4.5.　　The split-data approach, used for estimation of the nuisance functions and the error density, can be replaced by using all $n$ observations for estimation of all functions involved. The corresponding proof becomes much more complicated, and the interested reader is referred to [16].

4.6.　　All assertions hold if, in truncated cosine series estimates of the design density, regression and scale, defined in Section 3.4, the cutoff $S$ is changed on $n^{1/3}\ln(b_n)$. Then, under Assumption A, all these estimates are undersmoothed; that is, they have a bias which is smaller than an optimal one. This is an interesting remark for the reader who would like to understand the variance-bias balance in these estimates. Also, Efromovich [16] shows that rate-optimal adaptive estimation of nuisance functions can be also used. Thus, there is a robust choice among Fourier series estimates of the nuisance functions. On the other hand, it is an open problem to explore nonseries estimates like kernel or spline ones. Some numerical results can be found in [16].

4.7.　　For density estimation based on direct observations, there is a vast literature on closely related topics like censored observations, biased data, observations contaminated by measurement errors, estimation of functionals, ill-posed settings, estimation under a shape restriction, and so on. The obtained results indicate that it is reasonable to conjecture that many of the known direct-data results can be extended to the error density estimation as well. For instance, Van Keilegom and Veraverbeke [38] considered the problem of consistent error density estimation in a censored regression; using [4, 15, 18], it is reasonable to conjecture that optimal nonparametric results can be obtained for censored and biased regression models.

4.8. The reason for considering the case of dependent errors and predictors is threefold. First, this is a rather typical case in applications; the obtained result shows that in this case the marginal error density is exhibited by residuals. Second, we can conclude that the plug-in EP estimation is robust. Finally, let us stress that small datasets may not allow the statistician to evaluate a conditional density; then the univariate marginal density becomes a valuable tool for data analysis. Let us finish this remark by answering a question that the author was frequently asked during presentation of the result. Is it possible that the marginal error density is normal and the conditional density (of regression error given predictor) is not? The answer is "yes." As an example, define a bivariate density $\psi(u,x) := f(u) + \delta\lambda(u)\mu(x)$, where $f(u)$ is the standard normal density, $|\lambda(u)\mu(x)| \leq 1$, $\int_0^1 \mu(x)\,dx = \int_{-\infty}^{\infty} \lambda(u)\,du = 0$, and $\lambda(u) = 0$ whenever $f(u) < \delta$. There are plenty of such functions and, under the given assumptions, $\psi(u,x)$ is a valid bivariate density on $(-\infty,\infty) \times [0,1]$ with the standard normal marginal density $f(u)$. Obviously, the conditional density $\psi(u|x) := \psi(u,x)$ is not necessarily normal, and this verifies the assertion. The conclusion is that, even if the marginal distribution of residuals may be considered normal, unfortunately this does not imply the normality of the conditional distribution.

4.9. Brown, Low and Zhao [5] introduced the notion of nonparametric superefficiency, and they noticed that the Pinsker oracle (EP estimate) was superefficient; see also [18]. This fact, together with Theorem 1, immediately implies that the plug-in Pinsker oracle is also superefficient.

4.10. Let us note that plug-in estimation obviously enjoys its renascence in nonparametric estimation theory; see the discussion in [3] and [23]. A typical nonparametric plug-in setting is about optimal estimation of a functional. In this article a plug-in approach is caused by the indirect nature of observations, and, thus, it presents a new chapter in the theory of plug-in estimation.

4.11. It is a very interesting and technically involved problem to estimate the error density for the model with measurement errors in the predictors; see [7].

4.12.   The results hold for the case of a fixed-design regression; see [16].

4.13.   Let us comment on our main assumption about independence of pairs of observations $(X_1, Y_1), \ldots, (X_n, Y_n)$ with the typical counterexample being the case of dependent errors. The author conjectures that, based on the result of Efromovich [14], Section 4.8, even errors with a long memory should not affect the corresponding optimal rates. On the other hand, the outcome should change dramatically if a fixed design regression (say a time series) is considered. For this setting, the result of Hall and Hart [26] may be instrumental.

4.14.   It is an open and practically interesting topic to develop optimal Bayes and conditional distribution methods and then compare them with the developed plug-in estimator for the case of small datasets.

4.15.   Wavelet regression is a popular tool for solving many practical problems involving spatially inhomogeneous regressions. It is an open and interesting problem to explore the possibility of using wavelet-residuals as proxies for underlying regression errors.

## APPENDIX A: PROOFS

PROOF OF THEOREM 1.   Only the main steps of the proof are presented; the interested reader can find a detailed proof in [17, 19]. We begin with a more complicated case of finite support. Recall that the Pinsker oracle $\hat{f}_P$ is defined in Appendix B and it is based on pseudo-statistics $\{\bar{\mu}_k, \bar{\theta}_j\}$; in what follows we use the diacritics "bar" or "hat" above $\mu_k$ and $\theta_j$ to indicate a pseudo-statistic (oracle) based on underlying errors or a statistic based on observations, respectively. Set $\mathbf{Z}^* := (\varepsilon_{n-n_2+1}, \ldots, \varepsilon_n)$. Then a straight-forward calculation, based on $n_2 \geq [1 - 3(b_n^{-1} + n^{-1})]n$, establishes a plain inequality $E \int (\hat{f}_P(u, \mathbf{Z}^*) - f(u))^2 \, du \leq (1 + Cb_n^{-1}) E \int (\hat{f}_P(u, \mathbf{Z}) - f(u))^2 \, du$. As a result, in what follows we are assuming that pseudo-statistics $\bar{\mu}_k$ and $\bar{\theta}_j$ are based on $\mathbf{Z}^*$ in place of $\mathbf{Z}$, that is, plugged-in residuals correspond to errors used by the Pinsker oracle. Also recall that the oracle uses EP blockwise-shrinkage with $L_k = k^2$ and $t_k = \ln^{-2}(2+k)$. Keeping this in mind and using the Parseval identity, we write

$$E \int_0^1 (\hat{f}_P(u, \hat{\mathbf{Z}}) - f(u))^2 \, du$$
$$= E \sum_{k=1}^{K} \sum_{j \in B_k} (\hat{\mu}_k \hat{\theta}_j - \theta_j)^2 + \sum_{k>K} \sum_{j \in B_k} \theta_j^2$$

$$= E \sum_{k=1}^{K} \sum_{j \in B_k} [(\bar{\mu}_k \bar{\theta}_j - \theta_j) + \bar{\mu}_k(\hat{\theta}_j - \bar{\theta}_j) + (\hat{\mu}_k - \bar{\mu}_k)\hat{\theta}_j]^2$$

$$+ \sum_{k>K} \sum_{j \in B_k} \theta_j^2$$

$$\leq \left[(1 + \ln^{-1}(b_n)) E \sum_{k=1}^{K} \sum_{j \in B_k} (\bar{\mu}_k \bar{\theta}_j - \theta_j)^2 + \sum_{k>K} \sum_{j \in B_j} \theta_j^2\right]$$

$$+ 2(1 + \ln(b_n))$$

(A.1)
$$\times \left[\sum_{k=1}^{K} \sum_{j \in B_k} E \bar{\mu}_k^2 (\hat{\theta}_j - \bar{\theta}_j)^2 + \sum_{k=1}^{K} \sum_{j \in B_k} E(\hat{\mu}_k - \bar{\mu}_k)^2 \hat{\theta}_j^2\right]$$

$$= (1 + \ln^{-1}(b_n)) E \int_0^1 (\hat{f}_P(u, \mathbf{Z}^*) - f(u))^2 \, du$$

$$+ 2(1 + \ln(b_n))$$

$$\times \left[\sum_{k=1}^{K} \sum_{j \in B_k} E \bar{\mu}_k^2 (\hat{\theta}_j - \bar{\theta}_j)^2 + \sum_{k=1}^{K} E(\hat{\mu}_k - \bar{\mu}_k)^2 \sum_{j \in B_k} \hat{\theta}_j^2\right].$$

We need to evaluate the second term on the right-hand side of (A.1). Recall that $\xi = b\varepsilon + a$ and write

$$\hat{\theta}_j - \bar{\theta}_j$$

$$= n_2^{-1} \sum_{l=3n_1+1}^{n} [\varphi_j([Y_l - \hat{m}(X_l)]/b\hat{\sigma}(X_l) - a/b) - \varphi_j(\varepsilon_l)]$$

$$= n_2^{-1} \sum_{l=3n_1+1}^{n} [\varphi_j([m(X_l) + \sigma(X_l)\xi_l - \hat{m}(X_l)]/b\hat{\sigma}(X_l) - a/b)$$

(A.2)
$$- \varphi_j(\varepsilon_l)]$$

$$= n_2^{-1} \sum_{l=3n_1+1}^{n} \left[\varphi_j\left(\varepsilon_l + \frac{m(X_l) - \hat{m}(X_l)}{b\hat{\sigma}(X_l)} + \xi_l \frac{\sigma(X_l) - \hat{\sigma}(X_l)}{b\hat{\sigma}(X_l)}\right)\right.$$

$$\left. - \varphi_j(\varepsilon_l)\right].$$

Using the Taylor expansion for the cosine function, we can write

$$(\hat{\theta}_j - \bar{\theta}_j)^2 = \left\{n_2^{-1} \sum_{l=3n_1+1}^{n} [-\pi j H_l 2^{1/2} \sin(\pi j \varepsilon_l)\right.$$

$$\begin{aligned}
&\quad - (1/2)(\pi j)^2 H_l^2 2^{1/2} \cos(\pi j \varepsilon_l) \\
&\quad + (1/6)(\pi j)^3 H_l^3 2^{1/2} \sin(\pi j \varepsilon_l) \\
&\quad + (1/24)(\pi j)^4 H_l^4 2^{1/2} \cos(\pi j \varepsilon_l) \\
&\quad + (\pi j)^5 H_l^5 \bar{\nu}_l \Big] \Big\}^2
\end{aligned}$$

$$\begin{aligned}
(A.3) \quad \leq C \Bigg[ &j^2 n_2^{-2} \Big\{ \sum_{l=3n_1+1}^{n} H_l \sin(\pi j \varepsilon_l) \Big\}^2 \\
&+ j^4 n_2^{-2} \Big\{ \sum_{l=3n_1+1}^{n} H_l^2 \cos(\pi j \varepsilon_l) \Big\}^2 \\
&+ j^6 n_2^{-2} \Big\{ \sum_{l=3n_1+1}^{n} H_l^3 \sin(\pi j \varepsilon_l) \Big\}^2 \\
&+ j^8 n_2^{-2} \Big\{ \sum_{l=3n_1+1}^{n} H_l^4 \cos(\pi j \varepsilon_l) \Big\}^2 \\
&+ j^{10} n_2^{-2} \Big\{ \sum_{l=3n_1+1}^{n} |H_l|^5 \Big\}^2 \Bigg].
\end{aligned}$$

In the first equality we denoted by $\bar{\nu}_l$'s generic random variables satisfying $|\bar{\nu}_l| < 1$, and

$$(A.4) \quad H_l := \frac{m(X_l) - \hat{m}(X_l)}{b\hat{\sigma}(X_l)} + \xi_l \frac{\sigma(X_l) - \hat{\sigma}(X_l)}{b\hat{\sigma}(X_l)}.$$

As we see, the analysis of $(\hat{\theta}_j - \bar{\theta}_j)^2$ is converted into the analysis of the nonparametric regression and scale estimates. Evaluations are lengthy and technically involved (see them in [17]), and they imply

$$(A.5) \quad \sum_{k=1}^{K} \sum_{j \in B_k} E(\hat{\theta}_j - \bar{\theta}_j)^2 \leq C b_n n^{-1}.$$

Note that $\bar{\mu}_k^2 \leq 1$, so we have evaluated the first sum on the right-hand side of (A.1). Now let us consider the second sum. Write

$$\sum_{j \in B_k} (\hat{\mu}_k - \bar{\mu}_k)^2 \hat{\theta}_j^2$$

$$= L_k \Bigg[ \frac{\hat{\Theta}_k}{\hat{\Theta}_k + n_2^{-1}} - \frac{\bar{\Theta}_k}{\bar{\Theta}_k + n_2^{-1}} \Bigg]^2$$

$$\times (\hat{\Theta}_k + n_2^{-1}) I(\hat{\Theta}_k > t_k n_2^{-1}) I(\bar{\Theta}_k > t_k n_2^{-1})$$

(A.6)
$$+ \frac{\bar{\Theta}_k^2}{(\bar{\Theta}_k + n_2^{-1})^2} \sum_{j \in B_k} \hat{\theta}_j^2 I(\hat{\Theta}_k \leq t_k n_2^{-1}) I(\bar{\Theta}_k > t_k n_2^{-1})$$

$$+ \frac{L_k \hat{\Theta}_k^2}{\hat{\Theta}_k + n_2^{-1}} I(\hat{\Theta}_k > t_k n_2^{-1}) I(\bar{\Theta}_k \leq t_k n_2^{-1})$$

$$=: D_1(k) + D_2(k) + D_3(k).$$

Here we have used the notation $\bar{\Theta}_k := L_k^{-1} \sum_{j \in B_k} (\bar{\theta}_j^2 - n_2^{-1})$ and $\hat{\Theta}_k := L_k^{-1} \sum_{j \in B_k} (\hat{\theta}_j^2 - n_2^{-1})$.

Let us consider, in turn, these three terms, beginning with $D_1(k)$. Skipping the indicator functions, we are going to evaluate

$$D_1^*(k) := L_k \left[ \frac{\hat{\Theta}_k}{\hat{\Theta}_k + n_2^{-1}} - \frac{\bar{\Theta}_k}{\bar{\Theta}_k + n_2^{-1}} \right]^2 (\hat{\Theta}_k + n_2^{-1})$$

$$= \frac{L_k n_2^{-2} (\hat{\Theta}_k - \bar{\Theta}_k)^2}{(\hat{\Theta}_k + n_2^{-1})(\bar{\Theta}_k + n_2^{-1})^2}.$$

Using the Cauchy inequality, we can write, for any $c_k \geq 1$,

$$(\hat{\Theta}_k - \bar{\Theta}_k)^2 = L_k^{-2} \left[ \sum_{j \in B_k} (\hat{\theta}_j^2 - \bar{\theta}_j^2) \right]^2$$

(A.7)
$$\leq L_k^{-2} \left[ 2 c_k \sum_{j \in B_k} (\hat{\theta}_j - \bar{\theta}_j)^2 + c_k^{-1} \sum_{j \in B_k} \bar{\theta}_j^2 \right]^2$$

$$\leq 4 L_k^{-2} c_k^2 \left[ \sum_{j \in B_k} (\hat{\theta}_j - \bar{\theta}_j)^2 \right]^2 + 2 L_k^{-2} c_k^{-2} \left[ \sum_{j \in B_k} \bar{\theta}_j^2 \right]^2.$$

Note that $\sum_{j \in B_k} \bar{\theta}_j^2 = L_k(\bar{\Theta}_k + n_2^{-1})$, to get

$$D_1^*(k) \leq 4 n_2^{-2} \frac{(c_k^2 L_k^{-1})[\sum_{j \in B_k} (\hat{\theta}_j - \bar{\theta}_j)^2]^2}{(\hat{\Theta}_k + n_2^{-1})(\bar{\Theta}_k + n_2^{-1})^2} + 2 n_2^{-2} \frac{c_k^{-2} L_k}{\hat{\Theta}_k + n_2^{-1}}$$

$$=: D_{11}^*(k) + D_{12}^*(k).$$

Set $c_k^2 = L_k k^{1+d}$, $0 < d < 1$, and denote $D_{12}(k) := D_{12}^*(k) I(\hat{\Theta}_k > t_k n_2^{-1}) \times I(\bar{\Theta}_k > t_k n_2^{-1})$. We get

(A.8)
$$\sum_{k=1}^{K} D_{12}(k) \leq 2 n_2^{-1} \sum_{k=1}^{\infty} k^{-1-d} \frac{n_2^{-1} I(\hat{\Theta}_k > t_k n_2^{-1})}{\hat{\Theta}_k + n_2^{-1}} \leq C n^{-1}.$$

It is a more complicated task to evaluate $D_{11}^*(k)$. Denote

$$D_{11}(k) := D_{11}^*(k)I(\hat{\Theta}_k > t_k n_2^{-1})I(\bar{\Theta}_k > t_k n_2^{-1})$$

and write

(A.9) $$ED_{11}(k) \leq Cn_2 k^{1+d} E\left[\sum_{j \in B_k} (\hat{\theta}_j - \bar{\theta}_j)^2\right]^2.$$

The squared difference $(\hat{\theta}_j - \bar{\theta}_j)^2$ was considered in (A.3), and a calculation yields

(A.10) $$E\left[\sum_{j \in B_k} (\hat{\theta}_j - \bar{\theta}_j)^2\right]^2 \leq C[n^{-2}n^{-3/15}\ln^{10}(n) + n_1^{-2}k^{-6}].$$

Using this inequality in (A.9) implies $\sum_{k=1}^K ED_{11}(k) \leq Cn^{-1}b_n^2$. In its turn, this, together with (A.8), yields

$$E\sum_{k=1}^K D_1(k) \leq \sum_{k=1}^K E\{D_{11}(k) + D_{12}(k)\} \leq Cn^{-1}b_n^2.$$

Now we consider the second term $D_2(k)$ in (A.6). Write

$$D_2(k) = \frac{\bar{\Theta}_k^2}{(\bar{\Theta}_k + n_2^{-1})^2} L_k(\hat{\Theta}_k + n_2^{-1})I(\hat{\Theta}_k \leq t_k n_2^{-1})I(\bar{\Theta}_k > t_k n_2^{-1})$$

$$\leq Cn_2^{-1} L_k \frac{\bar{\Theta}_k^2}{(\bar{\Theta}_k + n_2^{-1})^2}[I(t_k n_2^{-1} < \bar{\Theta}_k \leq 2t_k n_2^{-1})$$

$$+ I(\bar{\Theta}_k > 2t_k n_2^{-1})I(\bar{\Theta}_k - \hat{\Theta}_k > \bar{\Theta}_k/2)].$$

Recall that $c_k^2 = L_k k^{1+d}$, $0 < d < 1$. Then using (A.7), we get

(A.11) $$(\hat{\Theta}_k - \bar{\Theta}_k)^2 \leq CL_k^{-2}c_k^2\left[\sum_{j \in B_k}(\hat{\theta}_j - \bar{\theta}_j)^2\right]^2 + Cc_k^{-2}(\bar{\Theta}_k + n_2^{-1})^2.$$

This, together with Chebyshev's inequality and (A.10), yields

$$ED_2(k) \leq Cn_2^{-1} L_k t_k^2 E\{I(t_k n_2^{-1} < \bar{\Theta}_k \leq 2t_k n_2^{-1})\}$$

$$+ Cn_2^{-1} L_k E\{I(\bar{\Theta}_k - \hat{\Theta}_k > \bar{\Theta}_k/2)I(\bar{\Theta}_k > 2t_k n_2^{-1})\}$$

(A.12) $$\leq Cn_2^{-1} L_k t_k^2 E\{I(t_k n_2^{-1} < \bar{\Theta}_k \leq 2t_k n_2^{-1})\}$$

$$+ Ck^{1+d}t_k^{-2}(n^{-1}n^{-3/15}\ln^{10}(n) + n_1^{-1}b_n k^{-6})$$

$$+ Cn_2^{-1}t_k^{-2}k^{-1-d}.$$

Let us evaluate the term $E\{I(t_k n_2^{-1} < \bar{\Theta}_k \leq 2t_k n_2^{-1})\}$. Denote $\Theta_k := L_k^{-1} \times \sum_{j \in B_k} \theta_j^2$ and $\theta_j := \int_0^1 f(u)\varphi_j(u)\,du$. Then using Lemma 1 in [11], together with some algebra, implies

$$E\{I(t_k n_2^{-1} < \bar{\Theta}_k \leq 2t_k n_2^{-1})\}$$
$$\leq I((1/2)t_k n_2^{-1} < \Theta_k \leq 4t_k n_2^{-1})$$
$$+ E\{I(\bar{\Theta}_k - \Theta_k > (1/2)t_k n_2^{-1})\}I(\Theta_k \leq (1/2)t_k n_2^{-1})$$
$$+ E\{I(\Theta_k - \bar{\Theta}_k > (1/2)\Theta_k)\}I(\Theta_k > 4t_k n_2^{-1})$$
$$\leq C t_k^{-1} \frac{\Theta_k}{\Theta + n_2^{-1}} I((1/2)t_k n_2^{-1} < \Theta_k \leq 4t_k n_2^{-1}) + C t_k^{-4} L_k^{-2}.$$

Let us recall a familiar blockwise Wiener oracle, which knows regression errors and an estimated density of errors and employs optimal shrinkage coefficients $\mu_k = \Theta_k/(\Theta_k + n^{-1})$. The Wiener oracle is the benchmark for the Pinsker oracle, and its MISE is proportional to $n_2^{-1} \sum_{k=1}^\infty L_k \Theta_k/(\Theta_k + n_2^{-1})$; see the discussion in [11, 14] and Appendix B. Then, combining the results, we get

$$\sum_{k=1}^K ED_2(k) \leq C n_2^{-1} \sum_{k=1}^K t_k L_k \frac{\Theta_k}{\Theta_k + n_2^{-1}} I((1/2)t_k n_2^{-1} < \Theta_k \leq 4t_k n_2^{-1}) + C b_n^2 n^{-1}$$
$$\leq C n_2^{-1} b_n^2 + C \sum_{k > b_n^{2/3}} t_k L_k \frac{\Theta_k}{\Theta_k + n_2^{-1}} I((1/2)t_k n_2^{-1} \Theta_k \leq 4t_k n_2^{-1})$$
$$+ C b_n^2 n^{-1}$$
$$\leq C \ln^{-2}(b_n) E \int_0^1 (\hat{f}_P(u, \mathbf{Z}^*) - f(u))^2 \, du + C b_n^2 n^{-1}.$$

Now we consider the third term $D_3(k)$ in (A.6). Write

$$D_3(k) \leq L_k \hat{\Theta}_k I(\hat{\Theta}_k > 2t_k n_2^{-1}) I(\bar{\Theta}_k \leq t_k n_2^{-1})$$
$$+ 2t_k L_k \hat{\Theta}_k I(t_k n_2^{-1} < \hat{\Theta}_k \leq 2t_k n_2^{-1}) I(\bar{\Theta}_k \leq t_k n_2^{-1})$$
$$=: D_{31}(k) + D_{32}(k).$$

To evaluate $D_{31}(k)$, we note that

$$D_{31}(k) \leq 2L_k |\hat{\Theta}_k - \bar{\Theta}_k| I(\hat{\Theta}_k - \bar{\Theta}_k > t_k n_2^{-1}) I(\bar{\Theta}_k \leq t_k n_2^{-1}).$$

This relation, the Chebyshev inequality, (A.10) and (A.11) imply

$$ED_{31}(k) \leq C L_k n_1^{-2} [L_k^{-2} c_k^2 (n^{-3/15} \ln^{10}(n) + k^{-6}) + c_k^{-2}]/(t_k n_2^{-1})$$
$$\leq C n_1^{-2} n_2 t_k^{-1} [k^{1+d}(n^{-3/5} \ln^{10}(n) + k^{-6}) + k^{-1-d}]$$

and, thus,

$$\sum_{k=1}^{K} ED_{31}(k) \leq Cn_1^{-2} n_2 < Cb_n^2 n^{-1}.$$

To evaluate $D_{32}(k)$, we write

$$\begin{aligned}
D_{32}(k) &= 2t_k L_k \hat{\Theta}_k I(t_k n_2^{-1} < \hat{\Theta}_k \leq 2t_k n_2^{-1}) I(\bar{\Theta}_k \leq (1/2)t_k n_2^{-1}) \\
&\quad + 2t_k L_k \hat{\Theta}_k I(t_k n_2^{-1} < \hat{\Theta}_k \leq 2t_k n_2^{-1}) I((1/2)t_k n_2^{-1} < \bar{\Theta}_k < t_k n_2^{-1}) \\
&=: D_{321}(k) + D_{322}(k).
\end{aligned}$$

The first term is evaluated in the same way as the term $D_{31}$ was evaluated, and we get $\sum_{k=1}^{K} ED_{321}(k) < Cb_n^2 n^{-1}$. The second term can be estimated as follows. First, we note that $\hat{\Theta}_k I(t_k n_2^{-1} < \hat{\Theta}_k \leq 2t_k n_2^{-1}) \leq 2t_k n_2^{-1}$. Then we realize that the term was evaluated earlier; see the first term in (A.12). This implies $\sum_{k=1}^{K} ED_{322}(k) \leq C \ln^{-2}(b_n) E \int_0^1 (\hat{f}_P(u, \mathbf{Z}^*) - f(u))^2 \, du + Cb_n^2 n^{-1}$. Combining the results, we get

$$\sum_{k=1}^{K} ED_3(k) \leq C \ln^{-2}(b_n) E \int_0^1 (\hat{f}_P(u, \mathbf{Z}^*) - f(u))^2 \, du + Cb_n^2 n^{-1}.$$

Then, by plugging the obtained estimates for $D_1(k)$, $D_2(k)$ and $D_3(k)$ into (A.6), we obtain

(A.13)
$$\begin{aligned}
E \sum_{k=1}^{K} \sum_{j \in B_k} (\hat{\mu}_j - \bar{\mu}_j)^2 \hat{\theta}_j^2 \\
\leq C \ln^{-2}(b_n) E \int_0^1 (\hat{f}_P(u, \mathbf{Z}^*) - f(u))^2 \, du + Cb_n^2 n^{-1}.
\end{aligned}$$

Using (A.5) and (A.13) in (A.1) verifies Theorem 1 for the finite support case. For the infinite support case, we set $\mathbf{Z}^* := (\varepsilon_{n-2n_1+1}, \ldots, \varepsilon_n)$, consider pseudo-statistics $\bar{h}(u)$ and $\bar{\mu}_k$ based on $\mathbf{Z}^*$, and then write, similarly to (A.1),

$$\begin{aligned}
&E \int_{-\infty}^{\infty} (\hat{f}_P(u, \hat{\mathbf{Z}}) - f(u))^2 \, du \\
&= E \int_{-\infty}^{\infty} \Bigg[ \pi^{-1} \sum_{k=1}^{K} \hat{\mu}_k \int_{B_k} \operatorname{Re}(\hat{h}(v) e^{-ivu}) \, dv \\
&\qquad\qquad - \pi^{-1} \sum_{k=1}^{K} \int_{B_k} \operatorname{Re}(h(v) e^{-ivu}) \, dv \\
&\qquad\qquad - \pi^{-1} \sum_{k>K} \int_{B_k} \operatorname{Re}(h(v) e^{-ivu}) \, dv \Bigg]^2 \, du
\end{aligned}$$

(A.14)
$$= \pi^{-1} E \sum_{k=1}^{K} \int_{B_k} |\hat{\mu}_k \hat{h}(v) - h(v)|^2 \, dv + \pi^{-1} \sum_{k>K} \int_{B_k} |h(v)|^2 \, dv$$
$$\leq (1 + \ln^{-1}(b_n)) E \int_{-\infty}^{\infty} (\hat{f}_P(u, \mathbf{Z}^*) - f(u))^2 \, du$$
$$+ 2(1 + \ln(b_n)) \bigg[ \sum_{k=1}^{K} E \bar{\mu}_k^2 \int_{B_k} |\hat{h}(v) - \bar{h}(v)|^2 \, dv$$
$$+ \sum_{k=1}^{K} E(\hat{\mu}_k - \bar{\mu}_k)^2 \int_{B_k} |\hat{h}(v)|^2 \, dv \bigg].$$

If we compare the right-hand side of (A.14) with the right-hand side of (A.1) and recall the steps taken after (A.1), then it is easy to recognize that the difference is in analyzing $\hat{h}(v) - \bar{h}(v)$ in place of $\hat{\theta}_j - \bar{\theta}_j$. Another remark is that now the second term in (A.4) vanishes because the unit scale function is known. These remarks allow us to follow the above-outlined proof and verify (3.6); details can be found in [19]. Theorem 1 is verified. □

PROOF OF COROLLARIES 1 AND 2. First, it is checked that the constant $P^*$ in (3.6) is uniformly bounded over the considered function classes. Second, it is easy to check that (3.6) holds with the estimate and the oracle exchanging places. Then using Theorem 1, Corollaries B1 and B2 of Appendix B, together with some algebra, verifies these corollaries. Details can be found in the technical reports. □

## APPENDIX B: THE PINSKER ORACLE

The Pinsker oracle is a data-driven density estimator possessing some desired statistical properties for the case of directly observed regression errors; in other words, it is a traditional density estimator whose oracle feature is in the knowledge of regression errors that are obviously unavailable to the statistician. For the case of direct observations and finite support, a good candidate for an estimator is the Efromovich–Pinsker (EP) data-driven (adaptive) procedure, which possesses an impressive bouquet of asymptotic properties of being: (a) minimax over a vast set of function classes which includes parametric, differentiable and analytic ones; (b) superefficient; (c) an excellent plug-in estimate; (d) applicable to filtering, regression and spectral density settings due to equivalence results. The interested reader can find discussion in [5, 9, 11, 14, 34, 39]. On the other hand, no results are available about a similar estimator for the case of a density with infinite support. The primary aim of this appendix is to develop such an estimator and explore its properties, and the secondary aim is to remind the reader of known results for the case of finite support.

We begin with the primary aim. Consider a density $f(z)$, $-\infty < z < \infty$, such that $\int_{-\infty}^{\infty} f^2(z)\,dx < \infty$. The problem is to estimate $f(z)$ under the MISE criteria when $n$ i.i.d. realizations $Z_1,\ldots,Z_n$ from $f$ are given. The underlying idea of EP estimation, translated from a finite support setting into the infinite one, is as follows. First, the characteristic function $h(v) := \int_{-\infty}^{\infty} e^{ivz} f(z)\,dz$ is estimated by its empirical counterpart $\hat{h}(v) := n^{-1}\sum_{l=1}^{n} e^{ivZ_l}$. Second, the estimate is "smoothed" by a statistic (filter) $\hat{\mu}(v)$, which is the main "ingredient" of the EP method defined shortly. Finally, a smoothed empirical characteristic function $\hat{\mu}(v)\hat{h}(v)$ is inverted to obtain a density estimate $\hat{f}(z) := (2\pi)^{-1}\int_{-\infty}^{\infty} \hat{\mu}(v)\hat{h}(v)e^{-ivz}\,dv$. Now we are in position to explain the underlying idea of choosing the EP smoothing. Consider a real even function $\mu(v):(-\infty,\infty)\to[0,1]$, set $\tilde{f}(z) := (2\pi)^{-1}\times\int_{-\infty}^{\infty}\mu(v)\hat{h}(v)e^{-ivz}\,dv$, and evaluate the MISE of this estimate using the Plancherel identity,

$$
\begin{aligned}
E\int_{-\infty}^{\infty} |\tilde{f}(z) - f(z)|^2\,dz & \\
(\text{B.1}) \quad &= (2\pi)^{-1} E\int_{-\infty}^{\infty} |\mu(v)\hat{h}(v) - h(v)|^2\,dv \\
&= (2\pi)^{-1} E\int_{-\infty}^{\infty} |\mu(v)(\hat{h}(v) - h(v)) - (1-\mu(v))h(v)|^2\,dv.
\end{aligned}
$$

Recall two familiar properties of the empirical characteristic function:

$$(\text{B.2}) \quad E\hat{h}(v) = h(v), \qquad E|\hat{h}(v) - h(v)|^2 = n^{-1}(1 - |h(v)|^2).$$

This, together with simple algebra, shows that a smoothing function (oracle),

$$(\text{B.3}) \quad \mu^*(v) := \frac{|h(v)|^2}{|h(v)|^2 + n^{-1}(1-|h(v)|^2)},$$

minimizes (B.1). The reader might notice that this smoothing function is the analog of the famous Wiener filter, and this is the reason why it also can be referred to as a filter; see [29], Chapter 10. The function $\mu^*(v)$ is unknown to the statistician, but, using (B.2), it can be estimated by the statistic

$$(\text{B.4}) \quad \tilde{\mu}(v) := \frac{|\hat{h}(v)|^2 - n^{-1}}{|\hat{h}(v)|^2} I(|\hat{h}(v)|^2 > (1+t)n^{-1}), \qquad t > 0;$$

here $I(\cdot)$ is the indicator function and $t$ is a threshold level ($1+t$ is often called a penalty). Hard thresholding (which is a trademark of the EP smoothing) is used to make the statistic a bona fide smoothing function. Unfortunately, it is not difficult to verify that this naive mimicry is not sufficiently accurate. Thus, by recalling that any characteristic function $h(v)$ is continuous and, thus, $\mu^*(v)$ is continuous, it is natural to approximate $\mu^*(v)$ by a piecewise constant function and then estimate that function. This is the underlying idea of the EP blockwise procedure. Note that $\mu^*(v)$ is an even function, and this allows us to work only with $v \in [0, \infty)$. We divide the half-line $[0, \infty)$ into a sequence of nonoverlapping blocks (intervals) $B_1, B_2, \ldots$ with corresponding lengths $L_k := \int_{B_k} dv > 0$, and then consider a smoothed empirical characteristic function $\tilde{h}(v) = \sum_{k=1}^{\infty} \mu_k \hat{h}(v) I(v \in B_k)$. Similarly to (B.1)–(B.3), we can establish that the MISE of the corresponding density estimate is minimized by the oracle

$$\text{(B.5)} \qquad \mu_k^* := \frac{L_k^{-1} \int_{B_k} |h(v)|^2 \, dv}{L_k^{-1} \int_{B_k} |h(v)|^2 \, dv + n^{-1}(1 - L_k^{-1} \int_{B_k} |h(v)|^2 \, dv)}.$$

Note the striking similarity between (B.3) and (B.5). Similarly to (B.4), the proposed estimate of the optimal $\mu_k^*$ is

$$\text{(B.6)} \qquad \hat{\mu}_k := \frac{L_k^{-1} \int_{B_k} |\hat{h}(v)|^2 \, dv - n^{-1}}{L_k^{-1} \int_{B_k} |\hat{h}(v)|^2 \, dv} I\left( L_k^{-1} \int_{B_k} |\hat{h}(v)|^2 \, dv > (1 + t_k) n^{-1} \right),$$

$$t_k > 0.$$

Then the EP density estimate is defined as

$$\text{(B.7)} \qquad \hat{f}^*(z) := \pi^{-1} \sum_{k=1}^{K} \hat{\mu}_k \int_{B_k} \mathrm{Re}(\hat{h}(v) e^{-ivz}) \, dv,$$

where the cutoff $K$ is a minimal integer such that $\sum_{k=1}^{K} L_k \geq n^{1/5} b_n$; this cutoff corresponds to the considered class of at least twice differentiable densities. The estimator (B.7) will be called the EP estimator for the case of infinite support. To better appreciate it, let us recall the EP density estimator for the finite support $[0, 1]$. The main difference is that here a discrete Fourier transform is used:

$$\text{(B.8)} \qquad \bar{f}^*(z) := 1 + \sum_{k=1}^{K} \bar{\mu}_k \sum_{j \in B_k} \bar{\theta}_j \varphi_j(z),$$

where $\{1, \varphi_j(z) = 2^{1/2} \cos(\pi j z), j = 1, 2, \ldots\}$ is the classical cosine basis on $[0, 1]$, $\{\bar{\theta}_j\}$ are empirical Fourier coefficients [estimates of Fourier coefficients $\theta_j := \int_0^1 f(z) \varphi_j(z) \, dz$]

$$\text{(B.9)} \qquad \bar{\theta}_j := n^{-1} \sum_{l=1}^{n} \varphi_j(Z_l),$$

and the smoothing weights (coefficients, filter) are

$$\text{(B.10)} \qquad \bar{\mu}_k := \frac{L_k^{-1} \sum_{j \in B_k} \bar{\theta}_j^2 - n^{-1}}{L_k^{-1} \sum_{j \in B_k} \bar{\theta}_j^2} I\left( L_k^{-1} \sum_{j \in B_k} \bar{\theta}_j^2 > (1 + t_k) n^{-1} \right).$$

Here the set of positive integers is divided into a sequence of blocks (including only neighbors) $B_k$ of cardinality $L_k$. Note that the EP infinite- and finite-support density estimates do look alike.

Finally, if $Z_1, \ldots, Z_n$ are unobserved regression errors, then the EP estimate becomes a *Pinsker oracle*. In this article, for both finite and infinite supports, this oracle is denoted as $\hat{f}_P(z, Z_1^n)$. Also, let us introduce the notation $\Theta_k := L_k^{-1} \int_{B_k} |h(v)|^2 \, dv$ and $\Theta_k := L_k^{-1} \sum_{j \in B_k} \theta_j^2$ for the infinite and finite supports, respectively, and $\mu_k := \Theta_k/(\Theta_k + n^{-1})$. Then $\hat{f}_P^*(z, Z^n)$ will denote a super-oracle (Wiener filter) which uses $\mu_k$ in place of $\hat{\mu}_k$ or $\bar{\mu}_k$ in the EP estimate; note that the super-oracle knows an estimated density $f$ and this is the oracle that is traditionally considered in the case of direct observations; see [14].

In what follows $C$'s denote generic positive constants and it is understood that the oracles vanish beyond the unit interval in the finite support case. Let us formulate a main property of the EP estimate; the result is new for the case of infinite support and it is due to Efromovich [11] for finite support.

THEOREM B1. *Suppose that $Z_1, \ldots, Z_n$ are i.i.d. realizations from a square integrable density $f$ with known support that can be either $[0,1]$ or $(-\infty, \infty)$. Consider the case of bounded thresholds $t_k < C$. Then the MISE of the Pinsker oracle (EP estimate) $\hat{f}_P(z, Z_1^n)$ satisfies the upper bound (oracle inequality)*

$$E \int_{-\infty}^{\infty} (\hat{f}_P(z, Z_1^n) - f(z))^2 \, dz$$

$$\leq \min \Bigg( E \int_{-\infty}^{\infty} (\hat{f}_P^*(z, Z_1^n) - f(z))^2 \, dz,$$

(B.11)
$$c^* \Bigg[ n^{-1} \sum_{k=1}^{K} L_k \mu_k + \sum_{k>K} L_k \Theta_k \Bigg] \Bigg)$$

$$+ \Bigg[ C n^{-1} \sum_{k=1}^{K} L_k \mu_k (t_k^{1/2} + L_k^{-1} t_k^{-3/2} d_f) \Bigg]$$

$$+ \Bigg[ C d_f^2 n^{-1} \sum_{k=1}^{K} L_k^{-1} t_k^{-3} \Bigg],$$

*where the constant $c^*$ is 1 or $\pi^{-1}$, and the functional $d_f$ is $1 + \sum_{j=1}^{\infty} |\theta_j|$ or $\int_0^{\infty} |h(v)| \, dv$, for the finite- and infinite-support cases, respectively.*

There are many important corollaries of Theorem B1. We present only two that are relevant to the topic of error density estimation: sharp minimax estimation of differentiable and analytic densities. We begin with the case of differentiable densities. For infinite support, we introduce a familiar Sobolev class $\mathcal{S}(\alpha, Q) := \{f : \int_{-\infty}^{\infty} [f^2(z) + (f^{(\alpha)}(z))^2] \, dz \leq Q\} = \{h : (2\pi)^{-1} \int_{-\infty}^{\infty} (1 + |v|^{2\alpha})|h(v)|^2 \, dv \leq Q\}$, where $f^{(\alpha)}$, $\alpha \geq 2$, is the $\alpha$th generalized derivative and $0 < Q < \infty$; see [33], page 144, [24] and [36]. With some obvious abuse of notation, for the case of finite support we define a similar Sobolev class $\mathcal{S}(\alpha, Q) := \{f : \sum_{j=1}^{\infty} (1 + (\pi j)^{2\alpha})\theta_j^2 \leq Q\}$; see [14], Chapter 2. Here we are interested only in the case $\alpha \geq 2$; more general Sobolev classes are considered in [20].

COROLLARY B1. *Consider the setting of Theorem* B1. *Suppose that $\alpha \geq 2$ and that blocks $B_k$ and thresholds $t_k$ used by the Pinsker oracle (EP estimate) $\hat{f}_P(z, Z_1^n)$ satisfy*

$$\text{(B.12)} \quad \sum_{k=1}^{\infty} L_k^{-1} t_k^{-3} < \infty, \qquad \lim_{k \to \infty} L_{k+1}/L_k = 1, \qquad \lim_{k \to \infty} t_k = 0.$$

*Then the Pinsker oracle (EP estimate) is sharp minimax over Sobolev densities and*

$$\text{(B.13)} \quad \begin{aligned} \sup_{f \in \mathcal{S}(\alpha, Q)} E \int_{-\infty}^{\infty} [r_n(\mathcal{S}(\alpha, Q))(\hat{f}_P(z, Z_1^n) - f(z))]^2 \, dz \\ = (1 + o(1)) \inf_{\check{f}} \sup_{f \in \mathcal{S}(\alpha, Q)} E \int_{-\infty}^{\infty} [r_n(\mathcal{S}(\alpha, Q))(\check{f}(z) - f(z))]^2 \, dz \\ = (1 + o(1)), \end{aligned}$$

*where the infimum is taken over all possible estimates $\check{f}$ based on observations $Z_1^n$ and parameters $\alpha$ and $Q$, and the sharp normalizing factor is defined in* (3.8).

Differentiable densities are traditionally studied in the nonparametric density estimation literature; see [24] and [37]. In the regression literature, typical error distributions are analytic, such as normal, mixture of normals and other stable distributions. For the case of infinite support, let us consider a class of such distributions studied in [25]. We say that $f$ belongs to an analytic class $\mathcal{A}(\gamma, Q)$, $0 < \gamma < \infty$, $0 < Q < \infty$, if $f(z)$, $-\infty < z < \infty$, has continuation into the strip $\{z + iy : |y| \leq \gamma, z \in (-\infty, \infty)\}$, $f(z + iy)$ is analytic inside this strip, bounded up to its boundary and $\int_{-\infty}^{\infty} (\text{Re}\{f(z + i\gamma)\})^2 \, dz \leq Q$. Note that this class includes, among others, normal, Student and Cauchy densities, as well as their mixtures and analytic one-to-one transformations. The main feature of these densities is a very fast (exponential) decrease of the corresponding characteristic functions, namely, according to Achieser ([1], page 251),

$$\text{(B.14)} \qquad \int_{-\infty}^{\infty} (e^{\gamma v} + e^{-\gamma v})^2 |h(v)|^2 \, dv \leq 8\pi Q.$$

As a result, we may expect almost parametric rates of MISE convergence. A finite-support counterpart of this class is well-known in the literature, and it is defined (with the obvious abuse of notation) as $\mathcal{A}(\gamma, Q) := \{f : \sum_{j=1}^{\infty} (1 + e^{2\pi \gamma j}) \theta_j^2 \leq Q\}$; see [14], Chapter 2.

COROLLARY B2. *Consider the setting of Theorem* B1, *and let* (B.12) *hold. Then the Pinsker oracle (EP estimate)* $\hat{f}_P(z, Z_1^n)$ *is sharp minimax over analytic densities and*

$$\text{(B.15)} \qquad \begin{aligned} &\sup_{f \in \mathcal{A}(\gamma, Q)} E \int_{-\infty}^{\infty} [r_n(\mathcal{A}(\gamma, Q))(\hat{f}_P(z, Z_1^n) - f(z))]^2 \, dz \\ &= (1 + o(1)) \inf_{\check{f}} \sup_{f \in \mathcal{A}(\gamma, Q)} E \int_{-\infty}^{\infty} [r_n(\mathcal{A}(\gamma, Q))(\check{f}(z) - f(z))]^2 \, dz \\ &= (1 + o(1)), \end{aligned}$$

*where the infimum is taken over all estimates* $\check{f}$ *based on observations* $Z_1^n$ *and parameters* $\gamma$ *and* $Q$, *and the sharp normalizing factor is defined in* (3.11).

These results show that the EP-estimate is simultaneously sharp adaptive over the union of differentiable and analytic densities. This allows us to conclude that the EP estimate is a feasible choice for a Pinsker oracle. Only to be specific, in Section 3 a Pinsker oracle with $L_k = k^2$ and $t_k = \ln^{-2}(2+k)$ is considered; note that this choice satisfies (B.12) and $t_k \leq 1$.

Now let us verify the stated results.

PROOF OF THEOREM B1. The assertion plainly follows from [11] for the finite-support case. Let us consider the infinite-support case. The plan is to follow along and employ the main parts of the proof presented in [11]; using the same notation will help us to do this. Set $\hat{\Theta}_k := L_k^{-1} \int_{B_k} |\hat{h}(v)|^2 \, dv - n^{-1}$ and note that (B.6) can be rewritten as $\hat{\mu}_k = \hat{\Theta}_k (\hat{\Theta}_k + n^{-1})^{-1} I(\hat{\Theta}_k > t_k n^{-1})$. Then, using (B.2) and the Plancherel identity, we write

$$\int_{-\infty}^{\infty} (\hat{f}_P(z, Z_1^n) - f(z))^2 \, dz$$

$$= \int_{-\infty}^{\infty} \left[ \pi^{-1} \sum_{k=1}^{K} \hat{\mu}_k \int_{B_k} \text{Re}(\hat{h}(v) e^{-ivz}) \, dv \right.$$

$$-\pi^{-1}\sum_{k=1}^{K}\int_{B_k}\mathrm{Re}(h(v)e^{-ivz})\,dv - \pi^{-1}\sum_{k>K}\int_{B_k}\mathrm{Re}(h(v)e^{-ivz})\,dv\Bigg]^2 dz$$

$$=\int_{-\infty}^{\infty}\left[\pi^{-1}\int_{0}^{\infty}\mathrm{Re}\left\{\left(\sum_{k=1}^{K}(\hat{\mu}_k\hat{h}(v)-h(v))I(v\in B_k)\right)e^{-ivz}\right\}dv\right]^2 dz$$

$$+\int_{-\infty}^{\infty}\left[\pi^{-1}\int_{0}^{\infty}\mathrm{Re}\left\{\sum_{k>K}h(v)I(v\in B_k)e^{-ivz}\right\}dv\right]^2 dz$$

$$=\pi^{-1}\sum_{k=1}^{K}\int_{B_k}|\hat{\mu}_k\hat{h}(v)-h(v)|^2\,dv + \pi^{-1}\sum_{k>K}\int_{B_k}|h(v)|^2\,dv.$$

This yields

$$E\int_{-\infty}^{\infty}(\hat{f}_P(z,Z_1^n)-f(z))^2\,dz$$

(B.16)
$$=\pi^{-1}\sum_{k=1}^{K}E\int_{B_k}|(\mu_k\hat{h}(v)-h(v))+(\hat{\mu}_k-\mu_k)\hat{h}(v)|^2\,dv$$
$$+\pi^{-1}\sum_{k>K}L_k\Theta_k$$
$$=:\pi^{-1}\sum_{k=1}^{K}A_k + \pi^{-1}\sum_{k>K}L_k\Theta_k.$$

Now we evaluate a particular $A_k$, $1\leq k\leq K$. Using the Cauchy inequality, we get

(B.17)
$$A_k \leq (1+t_k^{1/2})E\int_{B_k}|\mu_k\hat{h}(v)-h(v)|^2\,dv$$
$$+(1+t_k^{-1/2})E\left\{(\hat{\mu}_k-\mu_k)^2\int_{B_k}|\hat{h}(v)|^2\,dv\right\}$$
$$=:(1+t_k^{1/2})A_{k1} + (1+t_k^{-1/2})A_{k2}.$$

Note that $\pi^{-1}\sum_{k=1}^{K}A_{k1} + \pi^{-1}\sum_{k>K}L_k\Theta_k = E\int_{-\infty}^{\infty}(\hat{f}_P^*(z,Z_1^n)-f(z))^2\,dz$. On the other hand, using (B.2), we get

(B.18)
$$A_{k1} = E\int_{B_k}|\mu_k(\hat{h}(v)-h(v))-(1-\mu_k)h(v)|^2\,dv$$
$$= \mu_k^2 E\int_{B_k}|\hat{h}(v)-h(v)|^2\,dv + (1-\mu_k)^2\int_{B_k}|h(v)|^2\,dv$$

$$= n^{-1}\mu_k^2 \int_{B_k} (1 - |h(v)|^2)\, dv + (1 - \mu_k)^2 L_k \Theta_k$$

$$\leq L_k \left[ \frac{\Theta_k^2 n^{-1}}{(\Theta_k + n^{-1})^2} + \frac{n^{-2}\Theta_k}{(\Theta_k + n^{-1})^2} \right] = n^{-1} L_k \mu_k.$$

To evaluate $A_{k2}$, we note that $A_{k2} = L_k E\{(\hat{\mu}_k - \mu_k)^2 (\hat{\Theta}_k + n^{-1})\}$. Thus, at least formally, this term is identical to the same term in line (5.9) of [11]. To follow along the evaluation of $A_{k2}$ in [11], we need to verify that

$$E(\hat{\Theta}_k - \Theta_k)^4 \leq C \left( \int_{-\infty}^{\infty} |h(v)|\, dv \right)^2 L_k^{-2} n^{-2} (\Theta_k + n^{-1})^2,$$
(B.19)
$$1 \leq k \leq K.$$

This is done by a direct calculation which is similar to the proof of Lemma 3 in [11]; see also [20]. Then, similarly to lines (5.10)–(5.11) in [11], we get

$$\sum_{k=1}^{K} (1 + t_k^{-1/2}) A_{k2} \leq C n^{-1} \sum_{k=1}^{K} [L_k \mu_k (t_k^{1/2} + d_f L_k^{-1} t_k^{-3/2}) + d_f^2 L_k^{-1} t_k^{-3}].$$

Combining the results in (B.16) verifies Theorem B1. $\square$

PROOF OF COROLLARIES B1 AND B2. The second asymptotic equalities in these corollaries are established in [36]. The first asymptotic equalities follow from Theorem B1. $\square$

**Acknowledgments.** The author is grateful for the helpful and constructive comments of the Editor, Jianqing Fan, an Associate Editor and two referees.

DEPARTMENT OF MATHEMATICS AND STATISTICS
UNIVERSITY OF NEW MEXICO
ALBUQUERQUE, NEW MEXICO 87131
USA
E-MAIL: efrom@math.unm.edu